\documentclass[11pt,bezier]{article}
\usepackage{amsmath,graphicx,amssymb,amsfonts}
\renewcommand{\baselinestretch}{1.3}
\newtheorem{prethm}{{\bf Theorem}}

\newenvironment{thm}{\begin{prethm}{\hspace{-0.5
               em}{\bf.}}}{\end{prethm}}
\newtheorem{prepro}{{\bf Theorem}}

\newenvironment{pro}{\begin{prepro}{\hspace{-0.5
               em}{\bf.}}}{\end{prepro}}
\newtheorem{precor}{{\bf Corollary}}

\newtheorem{preconj}{{\bf Conjecture}}

\newtheorem{preremark}{{\bf Remark}}

\newtheorem{prelem}{{\bf Lemma}}

\newenvironment{lem}{\begin{prelem}{\hspace{-0.5
               em}{\bf.}}}{\end{prelem}}
\newtheorem{preproof}{{\bf Proof.}}

\newenvironment{proof}[1]{\begin{preproof}{\rm
               #1}\hfill{$\Box$}}{\end{preproof}}

\newcommand{\ch}{{\rm ch}}

\title{\bf\Large  {\sc Choice Number and Energy of Graphs}}

\author{{\normalsize{\sc Saieed Akbari}, {\sc Ebrahim
Ghorbani}}\vspace{2mm}
\\{\footnotesize{\it Department of Mathematical Sciences, Sharif
University of Technology,}}\vspace{-2mm}\\{\footnotesize{\it  P. O.
Box 11365-9415, Tehran, Iran}}\\{\footnotesize{\it Institute for
Studies in Theoretical Physics and
Mathematics,}}\vspace{-2mm}\\{\footnotesize{\it P. O. Box
19395-5746, Tehran, Iran}}\vspace{1mm}
\\{\footnotesize{$\mathsf{s\_akbari@sharif.edu}$\quad\quad
$\mathsf{e\_ghorbani@math.sharif.edu}$\vspace{-2mm}}}}

\date{}

\begin{document}
\maketitle

\begin{quote}
{\small \hfill{\rule{13.3cm}{.1mm}\hskip2cm}
\textbf{Abstract}\vspace{1mm}

{\renewcommand{\baselinestretch}{1}
\parskip = 0 mm

The energy of a graph $G$, denoted by $E(G)$, is defined as the sum
of the absolute values of all eigenvalues of  $G$. It is proved that
$E(G)\geq 2(n-\chi(\overline{G}))\ge 2(\ch(G)-1)$ for every graph
$G$ of order $n$, and that $E(G)\ge 2\ch(G)$ for all graphs $G$
except for those in a few specified families, where $\overline{G}$,
$\chi(G)$, and $\ch(G)$ are the complement, the chromatic number,
and the choice number of $G$, respectively. }}

\noindent{\small {\it Keywords}: Energy, choice number.}

\vspace{0.1mm}\noindent{\small {\it 2000 Mathematics Subject
Classification}:  05C15, 05C50, 15A03.

\vspace{-3mm}\hfill{\rule{13.3cm}{.1mm}\hskip2cm}}
\end{quote}

\vspace{9mm} \noindent{\bf\Large 1. Introduction}\vspace{3mm}

\noindent  All the graphs that we consider in this paper are finite,
simple and undirected. Let $G$ be  a graph. Throughout this paper
the {\it order} of $G$ is the number of vertices of $G$.  If $\{
v_1, \ldots, v_n\}$ is the set of vertices of $G$, then the {\it
adjacency matrix} of $G$, $A=[a_{ij}]$,
 is an $n\times{n}$ matrix where $a_{ij}=1$ if $v_i$ and $v_j$ are adjacent and $a_{ij}=0$
 otherwise. Thus $A$ is a symmetric matrix with zeros on the
 diagonal, and all the eigenvalues of $A$ are real and are denoted by
 $\lambda_1(G)\geq\cdots\geq\lambda_n(G)$. By the eigenvalues of $G$ we mean
 those of its adjacency matrix. The {\it energy} $E(G)$ of a graph $G$ is defined
 as the sum of the absolute
values of all eigenvalues of $G$, which is twice the sum of the
positive eigenvalues since the sum of all the eigenvalues is zero.
For a survey on the energy of graphs, see \cite{gut}.

 For a graph $G$, the {\it chromatic number} of $G$, denoted by
$\chi(G)$, is the minimum number of colors needed to color the
vertices of $G$ so that no two adjacent vertices have the same
color. Suppose that to each vertex $v$ of a graph G we assigned  a
set $L_v$ of $k$ distinct elements. If for any such assignment of
sets $L_v$ it is possible, for each $v\in V(G),$ to choose
$\ell_v\in L_v$ so that $\ell_u\neq\ell_v$ if $u$ and $v$ are
adjacent,  then $G$ is said to be {\it $k$-choosable}. The {\it
choice number} $\ch(G)$ of $G$ is the smallest $k$ such that $G$
is $k$-choosable.

We denote by $A_{n,t}$, $1\leq t \le n-1$, the graph obtained by
joining a new vertex to $t$ vertices of the complete graph $K_n$.
If we add two pendant vertices to a vertex of $K_n$, the resulting
graph has order $n+2$ and we denote it  by $B_n$.

In \cite{agz}, it is proved that apart from a few families of
graphs, $E(G)\geq 2\max(\chi(G), n-\chi(\overline{G}))$ (see the
following theorem). Our goal in this paper is to extend this result
to the choice number of graphs.

\begin{pro} \label{ab} Let $G$ be a graph. Then $E(G)<2\chi(G)$ if and
only if $G$ is a union of some isolated vertices and one of the
following graphs:\\ (i) the complete graph $K_n$;\\ (ii) the graph
$B_n$;\\
(iii) the graph $A_{n,t}$ for $n\leq7$, except when $(n,t)=(7,4)$,
and  also for $n\geq8$ and $t\in \{1, 2, n-1 \}$; \\ (iv) a triangle
with two pendant vertices adjacent to different vertices.
\end{pro}

The following is our main result.

\begin{thm} \label{main}  Let $G$ be a graph. Then $E(G)<2\,\ch(G)$ if and only if
$G$ is a union of some isolated vertices and one of the following
graphs:\\ (i)--(iv) as in Theorem \ref{ab};\\ (v) the complete
bipartite graph $K_{2,4}$.\end{thm}

\vspace{9mm} \noindent{\bf\Large 2. Proofs}\vspace{3mm}

In this section we present a proof for Theorem \ref{main}. To do so
we need some preliminaries.

 A well-known theorem of Nordhaus and Gaddum
\cite{nor} states that for every graph $G$ of order $n$,
$\chi(G)+\chi(\overline{G})\leq n+1$. This inequality can be
extended to the choice number. The graphs attaining equality are
characterized in \cite{dgm}. It is  proved that there are exactly
three types of such graphs defined as follows.
\begin{itemize}
  \item A graph  $G$ is of {\it   type $F_1$} if its vertex set can be partitioned
  into three sets $S_1, T, S_2$ (possibly, $S_2=\emptyset$) such that $S_1\cup S_2$ is an
  independent set of $G$, every vertex of $S_1$ is adjacent to every
  vertex of $T$, every vertex of $S_2$ has at least one non-neighbor
  in $T$, and $|S_1|$ is sufficiently large that
  the choice number of the induced subgraph on $T\cup S_1$ is
equal to $|T| + 1$. This implies that $\ch(G) = |T| + 1$ also.
Theorem 1 of \cite{gmm} states that if $T$ does not induce a
complete graph, then $|S_1|\ge |T|^2$; we will use this result
later.
  \item A graph is of {\it  type $\bar F_1$} if it is the complement of a graph of
   type $F_1$.
  \item A graph is of {\it  type $F_2$} if its vertex set can be
  partitioned into a clique $K$, an independent set $S$, and a 5-cycle $C$
  such that every vertex of $C$ is adjacent to every vertex of
  $K$ and to no vertex of $S$.

\end{itemize}

\begin{pro} \label{ch} (a) {\rm \cite{ert}} $\ch(G) + \ch(\overline{G})\le n + 1$ for every graph
$G$ of order $n$.\\ (b) {\rm \cite{dgm}} Equality holds in (a) if
and only if $G$ is of type $F_1$,  $\bar F_1$ or $F_2$.
\end{pro}

\begin{lem} \label{1}  For every graph $G$ of order $n$, $$E(G)\ge
2(n-\chi(\overline G))\ge2(n-\ch(\overline G))\ge2(\ch(G)-1).$$
\end{lem}

\begin{proof} {As remarked in \cite{agz}, the first inequality follows from
Theorem 2.30 of \cite{fav}, which states that
$n-\chi(\overline{G})\leq
\lambda_1(G)+\cdots+\lambda_{\chi(\overline{G})}(G)$.  The second
inequality holds because $\ch(G)\ge \chi(G)$ for every graph $G$,
and the third inequality holds by Theorem \ref{ch}$(a)$.}
\end{proof}




\begin{lem}\label{chrm} For every graph $G$, $\ch(G)\leq
\lambda_1(G)+1$.
\end{lem}

\begin{proof} {Wilf (\cite{w}, see also  \cite[p. 90]{spec}) proved that
every graph $G$ has a vertex with degree at most $\lambda_1(G)$, and
so does every induced subgraph of $G$. He deduced from this that
$\chi(G)\le \lambda_1(G)+1$, and the same  argument also proves that
$\ch(G)\le\lambda_1(G)+1$.}\end{proof}

\begin{lem} \label{2k2} Suppose $G$ has $2K_2$ as an induced subgraph. Then
$E(G)\ge2\ch(G)$. \end{lem}
\begin{proof} {By the Interlacing Theorem (Theorem
0.10 of \cite{spec}), $\lambda_2(G)\ge\lambda_2(2K_2)=1$, and so
$E(G)\ge2(\lambda_1(G)+\lambda_2(G))\ge2(\lambda_1(G)+1)\ge2\ch(G)$
by Lemma \ref{chrm}.}
\end{proof}

We are now in a position to prove Theorem \ref{main}.

\noindent\textbf{Proof of Theorem \ref{main}}. Let $G$ be a graph
such that $E(G) < 2\ch(G)$. We may assume that $G$ has at least one
edge, since otherwise $G$ is the union of some isolated vertices and
$K_1$, which is permitted by $(i)$ of Theorem \ref{main}. Since
removing isolated vertices does not change the value of $E(G)$ or
$\ch(G)$, we may assume that $G$ has no isolated vertices. If
$\ch(G) + \ch(\overline G) \le n$, then $E(G) \ge 2\ch(G)$ by Lemma
\ref{1}; this contradiction shows that $\ch(G) + \ch(\overline G) =
n + 1$, which means that $G$ has one of the types $F_1$, $\bar F_1$
and $F_2$ by Theorem \ref{ch}$(b)$. We consider these three cases
separately.


\noindent\textbf{Case 1}.  $G$ has type $F_1$. Then $G$ has
$G[T]\vee \overline{K}_k$ as an induced subgraph, where $G[T]$ is
the subgraph induced  by $G$ on $T$, $k=|S_1|$, and $\vee$ denotes
`join'. Let $|T|=t$, so that $\ch (G)=t+1$. If $G[T]$ is a complete
graph, then $\chi(G)=t+1=\ch(G)$, so that $E(G)<2\chi(G)$ and $G$ is
one of the graphs listed in Theorem \ref{ab}. So we may assume that
$G[T]$ is not a complete graph. In this case, as remarked after the
definition of type $F_1$, $k=|S_1|\ge|T|^2\geq t^2$. Thus
$$\lambda_1(G[T]\vee \overline{K}_k)\geq
\lambda_1(K_{t,t^2})=t\sqrt{t}\geq t+1,$$ provided $t\geq3$; since
$\ch(G)=t+1$, we have $E(G)\geq 2\,\ch(G)$. So we may assume that
$t\leq 2$, then $G[T]=\overline{K}_2$ and $k\geq t^2=4$.  For
$k\geq5$, we have $\lambda_1(K_{2,k})\geq \sqrt{10}>3=\ch(K_{2,k})$,
thus $E(G)\geq 2\,\ch(G)$. So we may assume that $k=4$. If $G\neq
K_{2,4}$, then either $|S_1|\ge5$ or $|S_2|>0$; thus $G$ has either
$K_{2,5}$ or $H$ as an induced subgraph, where $H$ is formed from
$K_{2,4}$ by adding an extra vertex joined to one of the vertices of
degree $4$. We have $E(K_{2,5})=2\sqrt{10}>6$. The graph $H$ has a
$P_4$ as an induced subgraph so $\lambda_2(H)\ge\lambda_2(P_4)>0.6$.
On the other hand $\lambda_1(H)\ge\lambda_1(K_{2,4})=2\sqrt{2}$.
Therefore $E(H)>2(2\sqrt{2}+0.6)>6$. Hence $E(G)>6=2\,\ch(G)$ if
$G\ne K_{2,4}$. Therefore $G=K_{2,4}$.

\noindent\textbf{Case 2}. $G$ has  type $\bar F_1$. So $\overline G$
is of  type $F_1$ with the associated partition $\{S_1, T, S_2\}$.
Let $t=|T|$ and $k=|S_1|$. If $\overline G[T]$ is not a complete
graph, then $k\ge t^2>1$ as in Case 1; hence $G$ has $2K_2$ as an
induced subgraph, which gives a contradiction by Lemma \ref{2k2}. So
$\overline G[T]$ is a complete graph.
Let $J$ be the set of those vertices of $T$ that are adjacent to all
vertices of $S_2$ in $G$. Let $v$ be a vertex of $S_1$. Then $G$ is
a graph of type $F_1$ with the
associated partition $\{S_1', T', S_2'\}$, in which 
$$
  \begin{array}{ll}
    S_1'=\{v\},~
T'=S_2\cup (S_1\setminus\{v\}),~ S_2'=T, & \hbox{if $k\ge2$;} \\
    S_1'=J\cup\{v\},~ T'=S_2,~ S_2'=T\setminus J, & \hbox{if $k=1$.}
  \end{array}$$
Therefore the result follows by Case 1.

\noindent\textbf{Case 3}. $G$ has  type $F_2$. Thus $G$ has a
5-cycle as an induced subgraph. So $\lambda_2(G)+\lambda_3(G)\geq
\lambda_2(C_5)+\lambda_3(C_5)>1$. 
Hence, by Lemma \ref{chrm},  we obtain
$$E(G)\geq 2(\lambda_1+\lambda_2+\lambda_3)>2(1+\lambda_1)\geq
2\ch(G).$$\hfill$\Box$



\noindent {\bf Acknowledgement.}  The  authors are indebted to the
Institute for Studies in Theoretical Physics and Mathematics (IPM)
for support; the research of the first author was in part supported
by a grant from IPM (No. 86050212). They are also grateful to the
referee for her/his helpful suggestions.

\end{document}